\documentclass[onefignum,onetabnum]{siamonline190516}

\usepackage{lipsum}
\usepackage{amsfonts}
\usepackage{graphicx}
\usepackage{epstopdf}
\usepackage{algorithmic}
\ifpdf
  \DeclareGraphicsExtensions{.eps,.pdf,.png,.jpg}
\else
  \DeclareGraphicsExtensions{.eps}
\fi

\usepackage{enumitem}
\setlist[enumerate]{leftmargin=.5in}
\setlist[itemize]{leftmargin=.5in}

\newsiamremark{remark}{Remark}
\newsiamremark{hypothesis}{Hypothesis}
\crefname{hypothesis}{Hypothesis}{Hypotheses}
\newsiamthm{claim}{Claim}

\headers{Dynamical organization of recollisions by a family of invariant tori}{J. Dubois, M. Jorba-Cusc\'{o}, \`{A}. Jorba, C. Chandre}

\title{Dynamical organization of recollisions by a family of invariant tori}

\author{J. Dubois \thanks{CNRS, Aix Marseille Univ, Centrale Marseille, I2M, Marseille, France. Current address: Max Planck Institute for the Physics of Complex Systems, N\"{o}thnitzer Stra\ss{e} 38, 01187 Dresden, Germany} \and M. Jorba-Cusc\'{o} \thanks{Departament de Matem\`{a}tiques i Inform\`{a}tica, Universitat de Barcelona, Barcelona, Spain. Current address: Universidad Internacional de la Rioja (UNIR), Logro\~{n}o, Spain} \and \`{A}. Jorba \thanks{Departament de Matem\`{a}tiques i Inform\`{a}tica, Universitat de Barcelona, Barcelona, Spain} \and C. Chandre \thanks{CNRS, Aix Marseille Univ, Centrale Marseille, I2M, Marseille, France}}

\usepackage{amsopn}

\begin{document}
\date{\today}
\maketitle


\begin{abstract}
We consider the motion of an electron in an atom subjected to a strong linearly polarized laser field. We identify the invariant structures organizing a very specific subset of trajectories, namely recollisions. Recollisions are trajectories which first escape the ionic core (i.e., ionize) and later return to this ionic core, for instance, to transfer the energy gained during the large excursion away from the core to bound electrons. We consider the role played by the directions transverse to the polarization direction in the recollision process. 
We compute the family of two-dimensional invariant tori associated with a specific hyperbolic-elliptic periodic orbit and their stable and unstable manifolds. We show that these manifolds organize recollisions in phase space.
\end{abstract}



\section{Introduction}
In strong-field atomic physics, recollisions are special trajectories of an electron of an atom or a molecule driven by an intense laser field in which the electron leaves the core by tunnel-ionization, gains energy from the laser field before coming back to the ionic core to collide with it~\cite{Corkum1993, Schafer1993}.
In essence, these trajectories are similar to ejection-collisions orbits in celestial mechanics~\cite{Olle2018_2}, in which the infinitesimal body (e.g., an asteroid) is first ejected from the primary body (e.g., a giant planet) before returning to collide with that primary body. In each of these types of trajectories, there is an ejection/ionization part and a subsequent collision part. 

Recollisions are the drivers of important nonlinear phenomena in atomic and molecular physics~\cite{Becker2008,Agostini2008,Krausz2009,Becker2012,Peng2015}, such as High Harmonic Generation (HHG)~\cite{Ferray1988,Lewenstein1994}, Above Threshold Ionization (ATI)~\cite{Paulus1994, Comtois2005} and Non-Sequential Multiple Ionization (NSMI)~\cite{Becker2012,Bergues2012}. These phenomena can be used to deduce information on the target atom or molecule, as a way to probe matter at the spatial and temporal scale of the electron motion~\cite{Corkum2007,Fohlisch2005,Goulielmakis2010,Haessler2010,Schultze2010,Klunder2011,Blaga2012,Peng2015}.
Here the questions we address are: How are these recollision trajectories organized in phase space~? In particular, what are the dynamical pathways or invariant structures guiding these trajectories away and back from the core~?  The answer to these questions holds the promise of controlling these strong-field phenomena by a more in-depth understanding of the influence of the laser parameters. 

A simplistic way to understand the origin of recollisions is to look at the motion of an electron of an atom along the polarization axis (with unit vector $\hat{\mathbf{x}}$) of a linearly polarized (LP) laser field, $\mathbf{E}(t) = \hat{\mathbf{x}} E_0 \cos ( \omega t)$. In the conventional recollision scenario~\cite{Corkum1993, Schafer1993, Krausz2009, Becker2008}, the electron first tunnel-ionizes through the Coulomb barrier induced by the laser at about the peak amplitude of the laser field. After tunnel-ionization, the electron moves classically in a combined Coulomb and laser field. 
Since the electron tunnels far away from the core, the electric force is dominant and the Coulomb force is neglected as a first approximation. The position of the electron along the polarization axis is approximated by
\begin{equation}
\label{eq:xsfa}
x(t) = x_0 + p_{x,0} t + ( E_0/\omega^2 ) \left[ \cos (\omega t ) - 1 \right].
\end{equation}
The resulting motion is a combination of a linear drift and an oscillatory motion due to the electric field. This oscillatory motion is responsible for pulling the electron away from the core and pulling it back half a laser cycle later, after an excursion of about $E_0/\omega^2$ away from the core. A rough condition for the return is that $p_{x,0}$ is not too large (typically smaller than $E_0/\omega$)~\cite{Bandrauk2005}.

If we now consider the directions transverse to the polarization axis, the electron motion is solely composed of a uniform drift, e.g., 
\begin{equation}
\label{eq:transverse_motion_SFA}
y(t) = y_0 + p_{y,0} t .
\end{equation}
For a non-zero initial momentum $p_{y,0}$, the electron drifts away from the core and never recollides in the strict sense. If $p_{y,0}$ is small enough (i.e., a small fraction of an atomic unit), the electron might return close enough to the core to trigger some meaningful energy exchange, provided that the position along the polarization axis is close to the core. 

This rough analysis suggests that recollisions originate from the dynamics along the polarization axis, and that the recollision dynamics along this axis is not structurally stable when considering the transverse directions, at least in absence of Coulomb interaction.
However, when the Coulomb interaction is fully taken into account, the electron dynamics in the transverse axis is coupled with its dynamics along the polarization axis. As a result, its dynamics along the transverse direction is significantly more complex than suggested by Eq.~\eqref{eq:transverse_motion_SFA} as we will show below. Here we study the nonlinear dynamics of the electron by identifying specific invariant objects as guides for recollisions. We focus on two families of periodic orbits: recolliding periodic orbits~\cite{Kamor2013,Kamor2014,Mauger2014_JPB,Abanador2017} (RPOs) and saddle periodic orbits~\cite{Mauger2012_PRL,Mauger2012_PRE,Barrabes2012,Olle2018,Olle2018_1} (SPOs) as relevant for the recollision process. We show that the RPOs lead to a fast escape in the transverse direction, i.e., they do not favor recollisions by themselves. However, a family of invariant tori associated with a specific SPO drives recollisions, and forms pathways for the motion of electrons away and back to the ionic core through their invariant manifolds. As a consequence, we show that the traditional picture obtained by neglecting the Coulomb field is highly inaccurate, especially when taking into account the dimension $d$ of the configuration space. The recollision scenario we detail below also demonstrates that the one-dimensional dynamics of the recollisions along the polarization axis is structurally stable when the Coulomb interaction is fully taken into account.

\section{Recollision scenario}

\subsection{The model}

We consider the classical motion of an electron in a combined electric and Coulomb field. For the Coulomb field, we consider an effective soft-Coulomb potential~\cite{Javanainen1988, Panfili2001, Panfili2002, Ho2005} to model the ground state of the electron. In atomic units, and in the dipole approximation, the Hamiltonian becomes
\begin{equation}
\label{eq:Hamiltonian}
    H({\bf r}, {\bf p},t) = \frac{{\bf p}^2}{2}-\frac{1}{\sqrt{{\bf r}^2+a^2}}+{\bf r}\cdot {\bf E}(t),
\end{equation}
where ${\bf r}, {\bf p} \in {\mathbb R}^d$. The above-Hamiltonian has $d+1/2$ degrees of freedom, so the only case whose dynamics can be easily represented (e.g., using Poincar\'e sections) is the one-dimensional case. It should be noted that the transverse subspaces $y=p_y=0$ and $z=p_z=0$ are invariant under the dynamics, which led to infer that the only relevant dynamics occurs along the polarization axis. The most interesting range of parameters occurs when the Coulomb interaction competes with the electric field, and none of the two can be considered as a perturbation of the other. In what follows, we consider an infrared laser field with 780 nm wavelength and an intensity of $3 \times 10^{14}$ W\ cm$^{-2}$. They correspond to $\omega=0.0584$ and $E_0 = 0.0925$ in atomic units (a.u.). The parameter $a$ measures how tightly bound an electron is in the ground state, or in other words, how easily this electron can be ionized by the field. We consider $a=1$  a.u.\ in what follows. 

Here we are investigating recollisions, i.e., trajectories which leave the vicinity of the ionic core located at ${\bf r}=0$, perform a large excursion away from it (at distances of the order of the quiver radius $E_0/\omega^2\approx 27$ a.u.\ or even further away) and later return to the vicinity of the ionic core (e.g., at distances of a few atomic units). 

In what follows, we first look at the one-dimensional case, and highlight a certain number of invariant structures which potentially play a role in the recollision process. From these structures, we elaborate possible scenarios driving recollisions. Then we look at which of these scenarios hold when considering the directions transverse to the polarization direction.  

The strategy here is to first look at the simplest invariant objects in phase space, namely periodic orbits. Given that the laser pulses are in general short (tens of laser cycles at most), we are looking at short periodic orbits. In addition, the interesting periodic orbits for recollisions should be weakly unstable: unstable because their stable/unstable manifolds guide the electron far away from and back to the core; weakly unstable because these orbits influence nearby motions for sufficiently long times to be most relevant.

\subsection{Recollisions in 1D}

The one-dimensional case ($d=1$) corresponds to $y=z=p_y=p_z=0$. Hamiltonian~\eqref{eq:Hamiltonian} becomes
\begin{equation}
\label{eq:Hamiltonian_1D}
H (x,p_x,t) = \dfrac{p_x^2}{2} - \dfrac{1}{\sqrt{x^2+1}} + x E_0 \cos (\omega t).
\end{equation}
\begin{figure}
\centering
\includegraphics[width=0.8\textwidth]{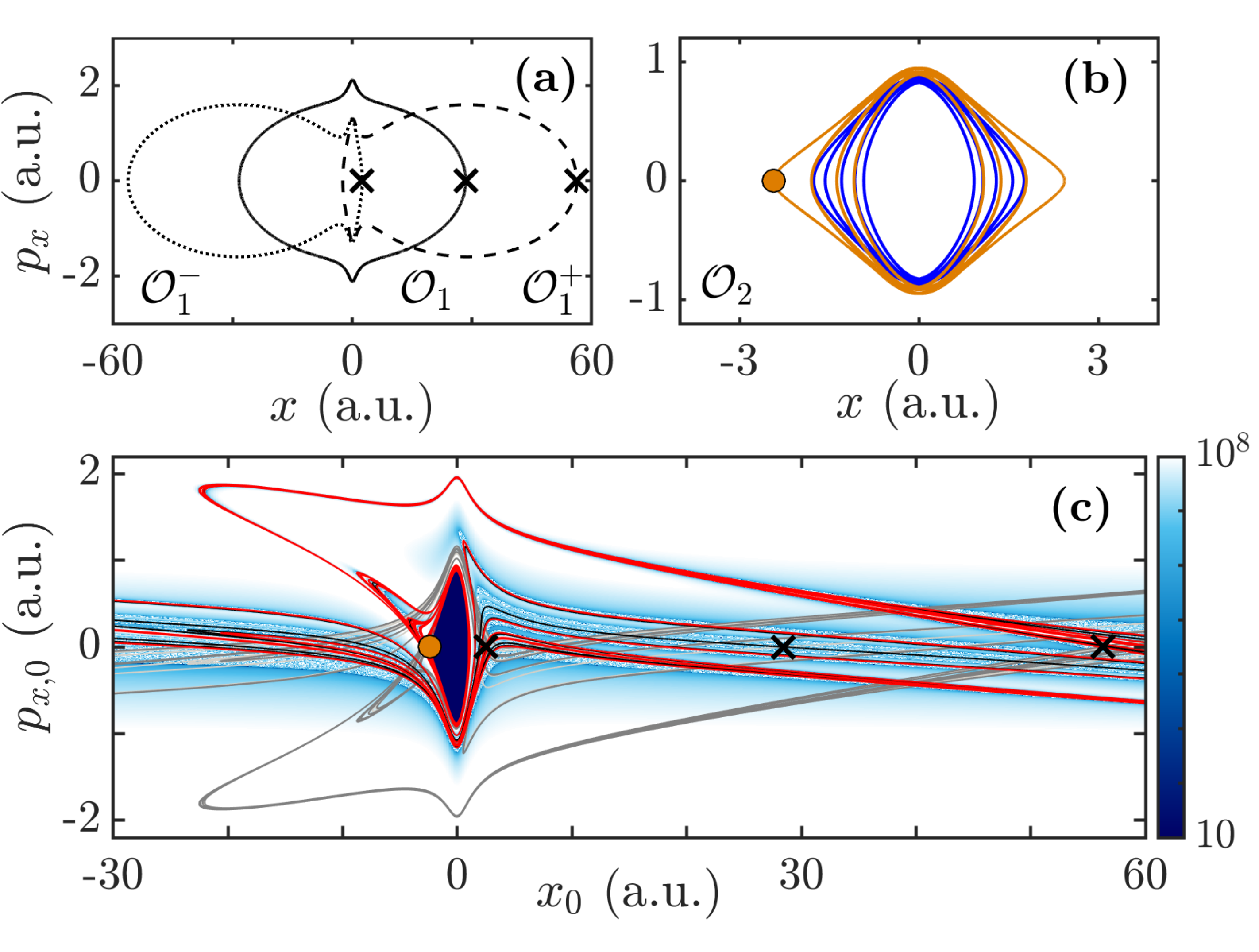}
\caption{Periodic orbit (a) $\mathcal{O}_1$, $\mathcal{O}_1^{\pm}$ and (b) $\mathcal{O}_2$ in the plane $(x,p_x)$ of Hamiltonian~\eqref{eq:Hamiltonian_1D} for the resonance $1:7$ (orange curve) and $1:9$ (blue curve). The crosses and the dot are the fixed points of the map ${\cal P}$ associated with these periodic orbits. (c) Distance between the electron and the ionic core, in logarithmic scale, for the trajectories of Hamiltonian~\eqref{eq:Hamiltonian_1D} as a function of the initial conditions $(x_0,p_{x,0})$ at $t=0$, for an integration time of $100 T$. The black and red lines are the stable manifold of $\mathcal{O}_1$ and $\mathcal{O}_2$, respectively. The light gray and gray lines are the unstable manifold of $\mathcal{O}_1$ and $\mathcal{O}_2$, respectively, which are the symmetric of their stable manifolds with respect to $p_x = 0$. The parameters are $E_0=0.0925$ and $\omega = 0.0584$. }
\label{fig:distance_manifold}
\end{figure}
\Cref{fig:distance_manifold}c shows the final distance (in logarithmic scale) of the electron to the ionic core as a function of the initial conditions $(x_0, p_{x,0})$ by integrating the Hamiltonian flow associated with Hamiltonian~\eqref{eq:Hamiltonian_1D} from $t = 0$ to $100 T$ where $T=2\pi/\omega$ is the laser period. We observe that the space of initial conditions is separated into two distinct regions: near the origin, a dark blue region corresponding to initial conditions for which the electron remains close to the ionic core, and a white/light blue region in which the electron eventually escapes/ionizes. In the latter region, we distinguish two different types of motion: A fast escape region which corresponds to a nearly-free motion [as in Eq.~\eqref{eq:xsfa}], located mainly in the white region, and a slow escape in which the motion is chaotic as shown by the sensitivity with respect to the initial conditions. The chaotic motion is a consequence of the recollisions and the nonlinearities experienced by the electron when it gets close to the ionic core. We observe that the regions of chaotic motion in the plane of initial conditions is organized in layers, and extend far away from the ionic core.
\par

\begin{figure}
\centering
\includegraphics[width=0.8\textwidth]{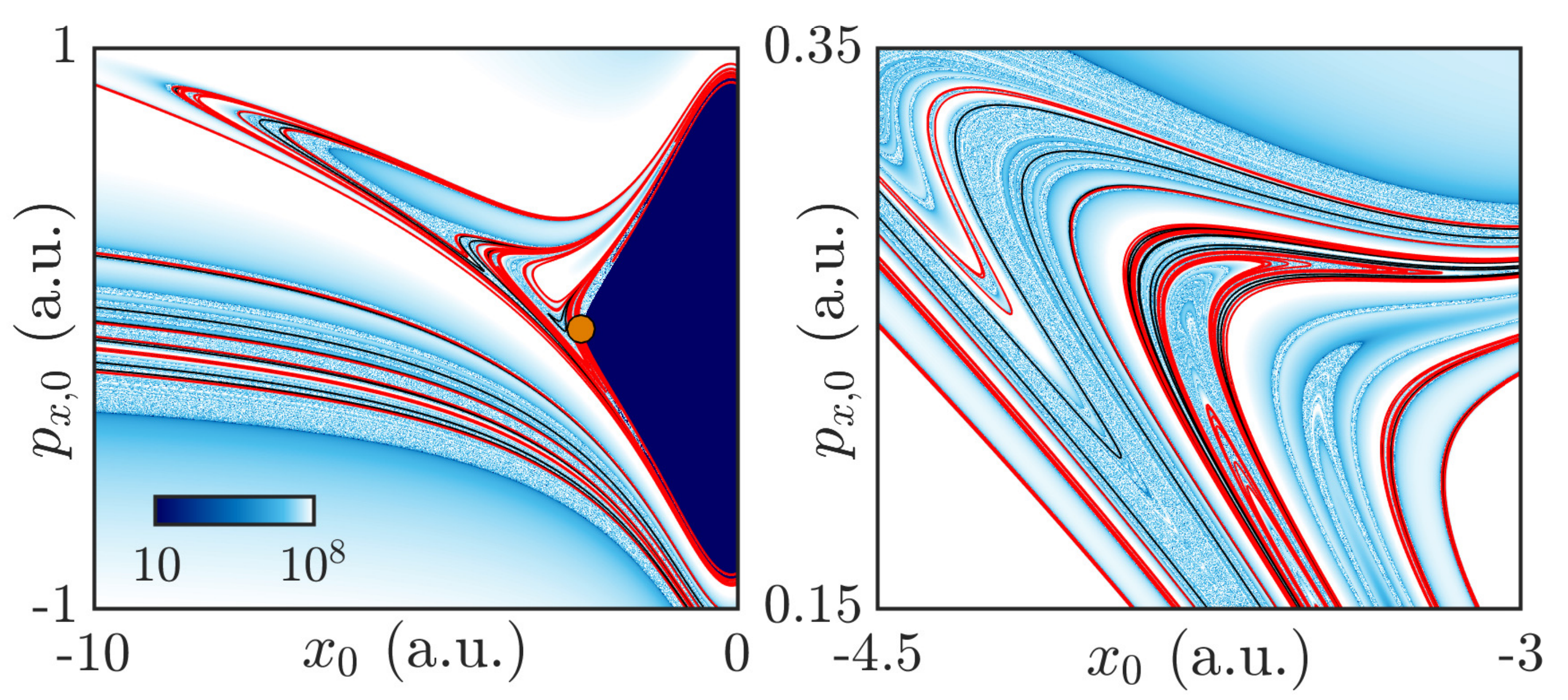}
\caption{Zoom of \cref{fig:distance_manifold}c close to the bounded region.}
\label{fig:distance_manifold_zoom}
\end{figure}

The basic organization of the dynamics can be understood by looking at weakly unstable periodic orbits of short duration, typically with period $T$. In what follows we mainly discuss four such periodic orbits, $\mathcal{O}_1$, $\mathcal{O}_1^{\pm}$ and $\mathcal{O}_2$. 
\begin{table}[]
    \centering
    \begin{tabular}{c|c c c}
         & $\mathcal{O}_1$ & $\mathcal{O}_1^{\pm}$ & $\mathcal{O}_2$ \\
         \hline\hline
        $\lambda_x^{(1)}$ & 2.3417 & $-16.410$  & 7.7463 \\
        $\lambda_x^{(2)}$ & 0.42705 & $-6.0937 \times 10^{-2}$  & 0.12909 \\
        $\lambda_{\rho}^{(1)}$ & $1.0662 \times 10^{3}$ & $1.1939 \times 10^{2}$ & $\exp (1.3499 \, {\rm i})$ \\
        $\lambda_{\rho}^{(2)}$ & $9.3788 \times 10^{-4}$ & $8.376\times 10^{-3}$ & $\exp (- 1.3499 \, {\rm i})$ \\
        $\lambda_{\theta}^{(1)}$ & 1 & 1 & 1 \\
        $\lambda_{\theta}^{(2)}$ & 1 & 1 & 1
    \end{tabular}
    \caption{Eigenvalues of the periodic orbits $\mathcal{O}_1$, $\mathcal{O}_1^{\pm}$ and $\mathcal{O}_2$ of resonance $1:7$ of Hamiltonian~\eqref{eq:Hamiltonian} for $a=1$, $E_0=0.0925$ and $\omega = 0.0584$.}
    \label{tab:eigenvalues_PO1D}
\end{table}
In Ref.~\cite{Kamor2014} it was argued that the invariant object organizing the recollisions in phase space is the unstable manifold associated with a RPO, which is here referred to as $\mathcal{O}_1$. Associated with $\mathcal{O}_1$ are the RPOs $\mathcal{O}_1^{\pm}$ which originate from a bifurcation of $\mathcal{O}_1$ as intensity increases. The RPOs can be seen as remainders of a family of periodic orbits in absence of Coulomb interaction. Indeed, by looking at Eq.~\eqref{eq:xsfa}, we notice that the motion is periodic for $p_{x,0}=0$, and for all $x_0$, in absence of Coulomb interaction. When the Coulomb interaction is increased (e.g., by increasing the effective charge), all but a finite number of these orbits are broken. We notice that some of these orbits are associated with an infinite number of recollisions since they periodically come back to the ionic core after large excursions far away from it. 
The RPOs $\mathcal{O}_1$ and $\mathcal{O}_1^{\pm}$ are depicted in black in \cref{fig:distance_manifold}a. These RPOs are weakly hyperbolic in the polarization plane since the monodromy matrix has one eigenvalue of modulus larger than one but still of order one (see \cref{tab:eigenvalues_PO1D}). An electron trajectory close to $\mathcal{O}_1$ exhibits the same pattern as a typical recolliding trajectory for a certain time: It crosses the origin $x=0$ when the amplitude of the electric field is at its lowest, and is far from the core at its peaks, twice a laser cycle. Under the Poincar\'e map ${\cal P} (x, p_x ) = \phi_0^{T} (x, p_x)$, where $\phi_0^{T}$ is the Hamiltonian flow from $t=0$ to $T$, the periodic orbit $\mathcal{O}_1$ is a fixed point of ${\cal P}$ with $x_0 \approx E_0/\omega^2$ and $p_{x,0}=0$. In \cref{fig:distance_manifold}c, we show the stable (black) and unstable (light gray) manifold of $\mathcal{O}_1$. The stable and unstable manifolds of the periodic orbits in 1D are computed using the Hobson method~\cite{Hobson1993}. The stable manifold of $\mathcal{O}_1$ reproduces very well the patterns of recollisions since we notice a very good agreement between this stable manifold and the ridges of the distance contour plot. The stable manifold is the one which brings the electron back close to the core, reducing the distance from the core even if the trajectories eventually ionize. It was then conjectured~\cite{Kamor2014} that the stable manifold of $\mathcal{O}_1$ is structuring the recollisions in phase-space. In particular, we see that the stable manifold of $\mathcal{O}_1$ crosses the $\hat{\mathbf{x}}$-axis where it coexists with the region of bounded motion. Notice that from the quantum mechanical point of view, this coexistence plays a crucial role, and finds its application in the HHG, since it gives rise to interferences between the bounded and recolliding parts of the wavepacket~\cite{Sand1999,Berman2018,Berman2019}.
\par
Regarding the SPOs, notably $\mathcal{O}_2$, they are the results of a resonant process $1:n$ with $n$ odd as described in Ref.~\cite{Mauger2012_PRE}. All have the same period, namely $T$. At the intensity we consider, the SPO associated with the resonance $1:7$ is the most relevant. 
This SPO $\mathcal{O}_2$ is depicted in \cref{fig:distance_manifold}b. A noticeable feature is that it is much closer to the core than $\mathcal{O}_1$ and $\mathcal{O}_1^{\pm}$. Strong evidence suggests that this periodic orbit plays an important role in the delayed double ionization~\cite{Mauger2012_PRL,Mauger2012_PRE}. It surrounds the bound region as observed in \cref{fig:distance_manifold}c, and its motion is not prototypical of recollisions in the sense that it does not present large excursions away from the core as compared with $ \mathcal{O}_1$ and $\mathcal{O}_1^{\pm}$. The periodic orbit $\mathcal{O}_2$ organizes the dynamics at the border of bounded and unbounded motion. As $\mathcal{O}_1$ and $\mathcal{O}_1^{\pm}$, the SPO $\mathcal{O}_2$ is weakly hyperbolic in the polarization plane (see \cref{tab:eigenvalues_PO1D}). In the vicinity of the region of bounded motion, the stable and unstable manifolds of $\mathcal{O}_2$ intersect each other, creating a thin chaotic tangle. Some electrons can escape through this tangle with some ionization delays. This is the mechanism for recollision-excitations with subsequent ionizations (RESI) for non-sequential multiple ionizations~\cite{Mauger2012_PRL,Mauger2012_PRE}. Relatively far away from the region of bounded motion, the stable manifold of $\mathcal{O}_2$ reproduces well the patterns of recollisions revealed by the contour plot of the final distance of the electron, as observed in \cref{fig:distance_manifold}c. In other words, the stable and unstable manifolds of $\mathcal{O}_1$ and $\mathcal{O}_2$ display similar behaviours, and therefore it can be equally argued that $\mathcal{O}_2$ also drives recollisions in 1D, even if $\mathcal{O}_2$ does not itself display recollision features. 
\par

For the one-dimensional dynamics, the stable and unstable manifolds of the RPOs $\mathcal{O}_1$, and $\mathcal{O}_1^{\pm}$ and the SPO $\mathcal{O}_2$ partition the phase space and drive the recollisions and the electron dynamics through their stable and unstable manifolds. Given that the one-dimensional dynamics corresponds to an invariant subspace, these periodic orbits persist in higher dimensions. Therefore they are good candidates for driving recollisions in higher dimensions $d=2,3$.
\begin{figure}
\centering
\includegraphics[width=0.8\textwidth]{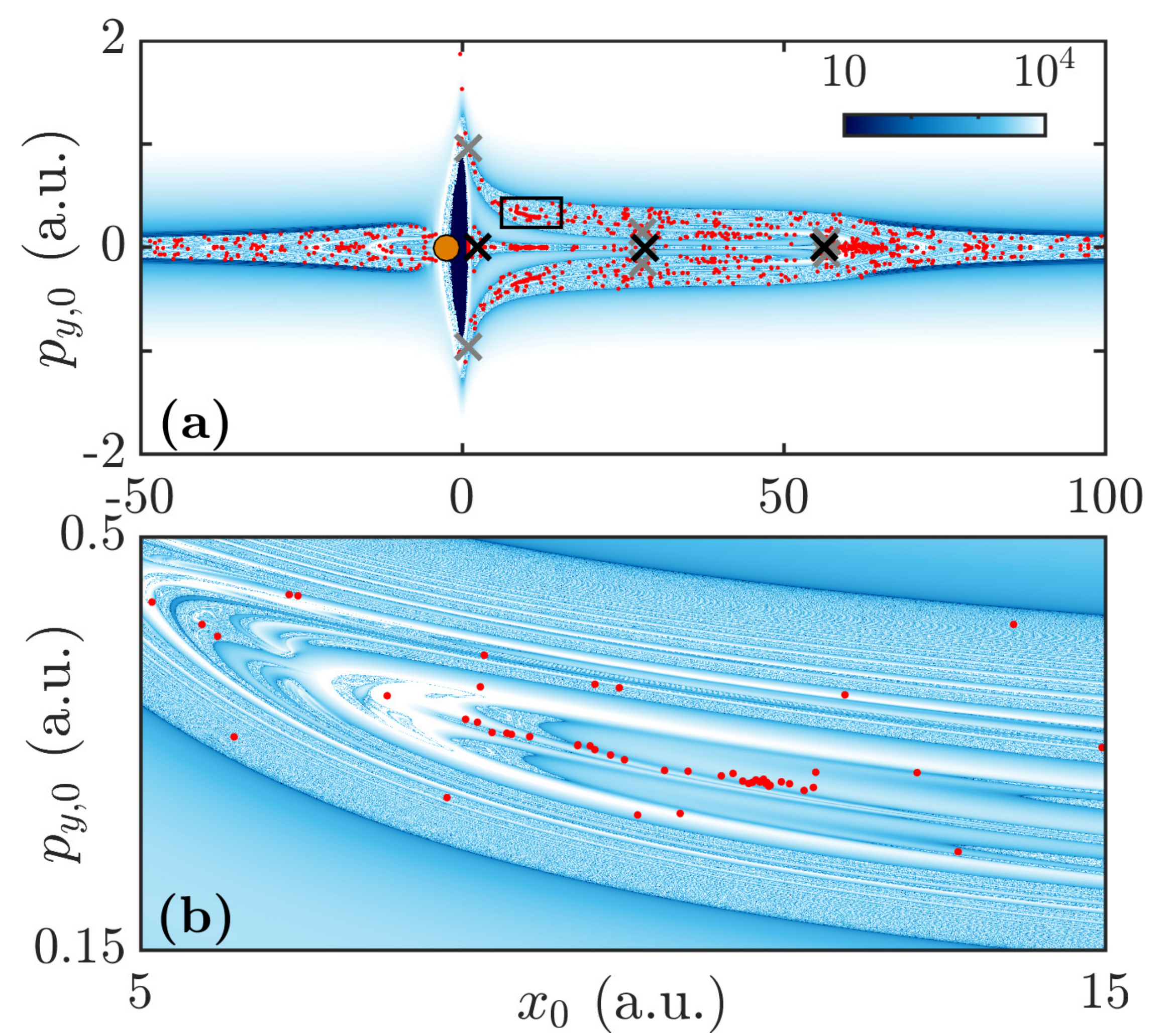}
\caption{Distance between the electron and the ionic core, in logarithmic scale, solution of Hamiltonian~\eqref{eq:Hamiltonian_2D} as a function of the initial conditions $(x_0,p_{y,0})$ and $y_0 = p_{x,0} = 0$ at $t=0$, for an integration $100 T$ and initial condition $t=0$. The red dots are points on the stable manifold of $\mathcal{F}$, the family of invariant tori associated with $\mathcal{O}_2$. The dark blue region are the initial conditions of an electron in a bound state. The black crosses, gray crosses and orange dot correspond to the fixed point associated with $\mathcal{O}_1$ and $\mathcal{O}_1^\pm$, $\mathcal{O}$ and $\mathcal{O}^\pm$, and $\mathcal{O}_2$ under $\mathcal{P}$, respectively. The panel (b) is a zoom of (a) corresponding to the black square in (a). The parameters are $E_0=0.0925$ and $\omega = 0.0584$.}
\label{fig:x_py_surface}
\end{figure}

\subsection{Recollisions in 2D}
\begin{figure}
\centering
\includegraphics[width=0.5\textwidth]{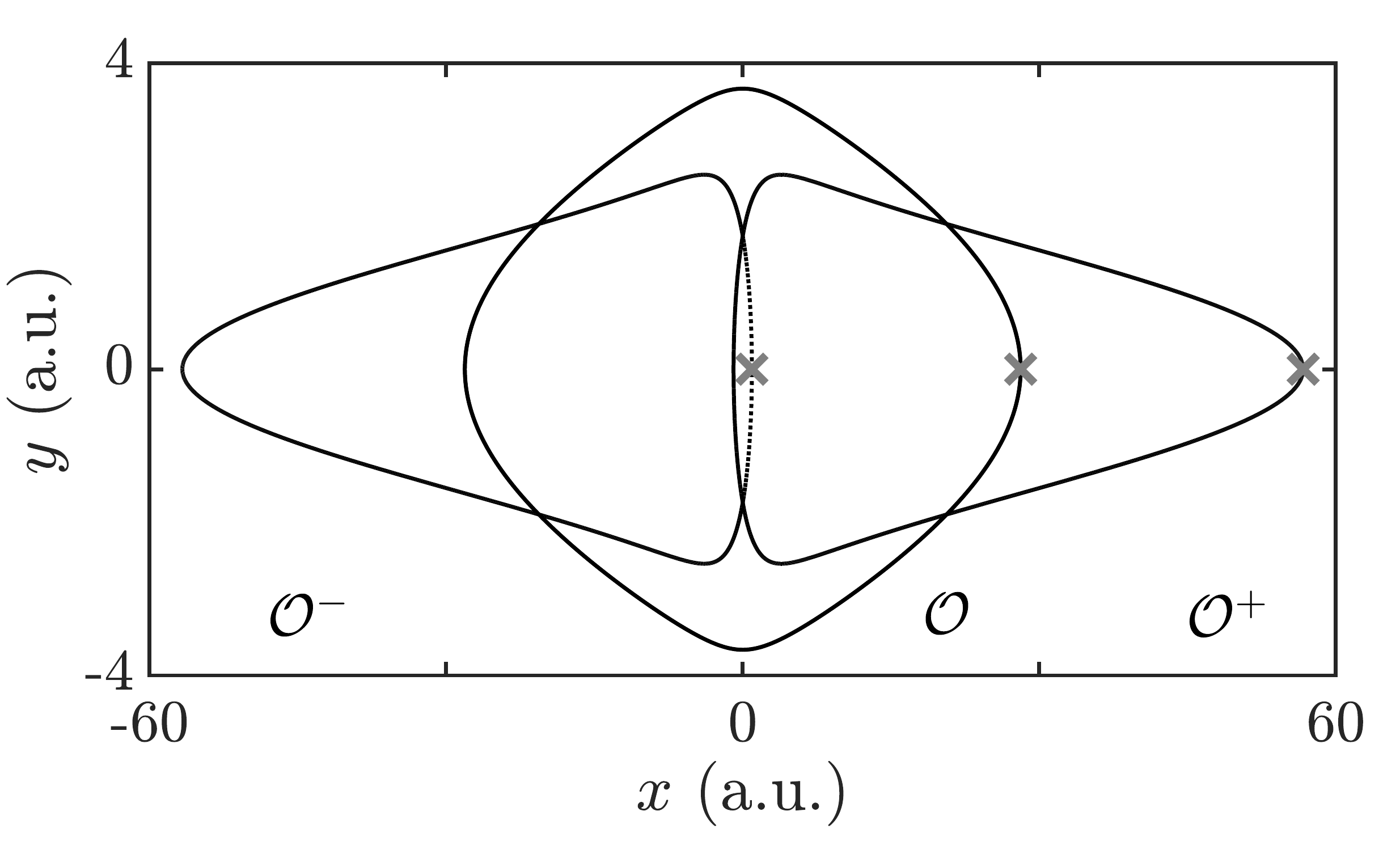}
\caption{Projections of the RPOs $\mathcal{O}$ and $\mathcal{O}^\pm$ of Hamiltonian~\eqref{eq:Hamiltonian_2D}. The crosses show the fixed points of ${\cal P}$. The parameters are $E_0=0.0925$ and $\omega = 0.0584$.}
\label{fig:RPO_2D}
\end{figure}
\begin{table}[]
    \centering
    \begin{tabular}{c|c c}
         & $\mathcal{O}$ & $\mathcal{O}^{\pm}$ \\
         \hline\hline
        $\lambda_1^{(1)}$ & 25.065 & $-29.910$ \\
        $\lambda_1^{(2)}$ & 2.1379 & $-11.438$ \\
        $\lambda_2^{(1)}$ & 0.46775 & $-0.033434$ \\
        $\lambda_2^{(2)}$ & 0.039896 & $-0.08743$ \\
        $\lambda_{\theta}^{(1)}$ & 1  & 1 \\
        $\lambda_{\theta}^{(2)}$ & 1  & 1
    \end{tabular}
    \caption{Eigenvalues of the RPOs $\mathcal{O}$ and $\mathcal{O}^{\pm}$ of Hamiltonian~\eqref{eq:Hamiltonian} for $a=1$, $E_0=0.0925$ and $\omega = 0.0584$.}
    \label{tab:eigenvalues_PO2D}
\end{table}

\subsubsection{Highly-dimensional objects drive the dynamics}

For $d=2$, Hamiltonian~\eqref{eq:Hamiltonian} becomes
\begin{equation}
\label{eq:Hamiltonian_2D}
H (x,y,p_x,p_y,t) = \dfrac{p_x^2}{2} + \dfrac{p_y^2}{2} - \dfrac{1}{\sqrt{x^2+y^2+1}} + x E_0 \cos (\omega t) .
\end{equation}
First, we observe that the degrees of freedom $(x,p_x)$ and $(y,p_y)$ are fully coupled through the ion-electron interaction term, which is especially strong close to the ionic core, but also significant far away due to the long-range interaction of the Coulomb potential. In \cref{fig:x_py_surface}, we show the final distance of the electron to the ionic core as a function of the initial conditions $(x_0, p_{y,0})$ with $y_0 = p_{x,0} = 0$. In a similar way as for $d=1$, we also observe two distinct regions in phase space: near the origin, a dark blue region corresponding to bounded motion, and a white/light blue region in which the electron motion is unbounded. In the latter region, there is an unbounded but mostly regular motion (white region), and a chaotic escape motion characterized by the sensitivity with respect to initial conditions. As for $d=1$, the chaotic region is organized in layers in the plane of initial conditions $(x_0,p_{y,0})$ and $y_0 = p_{x,0} = 0$ and is computed for a fixed initial time $t=0$, therefore it corresponds to a two-dimensional slice in a five-dimensional phase space. The fact that recollisions are organized in layers on this slice of initial conditions suggests that the domain of initial conditions leading to recollisions is linked with a codimension one manifold: From \cref{fig:x_py_surface}, it appears that the recollisions fill a nearly two-dimensional region in the space $(x_0, p_{y,0})$ around $p_{y,0}\approx 0$ which extends far away from the ionic core, the backbone of which is likely one-dimensional as for $d=1$. However, since the space of initial conditions $(x_0, y_0, p_{x,0}, p_{y,0})$ is four-dimensional, this implies that this backbone is three-dimensional. As a consequence, a good candidate for organizing recollisions for $d=2$ would be a three-dimensional invariant manifold for $\mathcal{P}(x,y,p_x,p_y)=\phi_0^T(x,y,p_x,p_y)$ where $\phi_0^T$ is the Hamiltonian flow from $t=0$ to $t=T$ associated with Hamiltonian~\eqref{eq:Hamiltonian_2D}.
\par

\begin{figure}
\centering
\includegraphics[width=\textwidth]{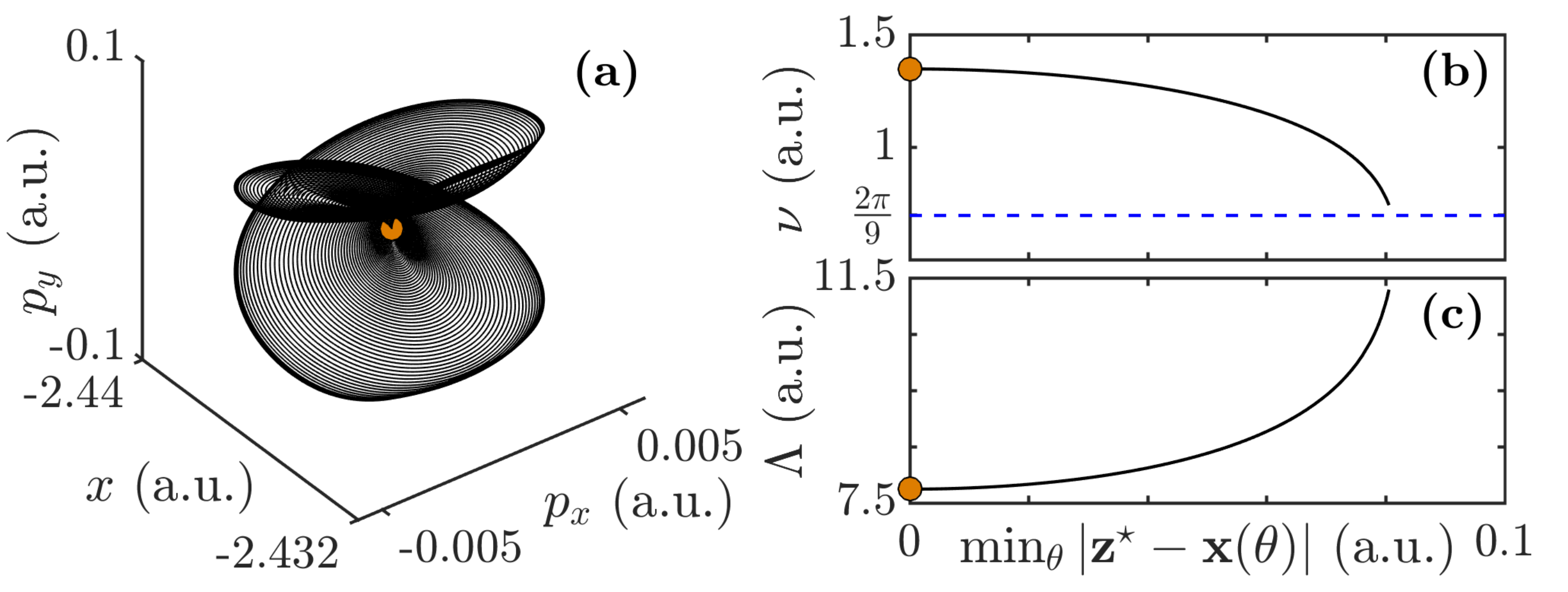}
\caption{(a) Invariant curves (black curves) of the fixed point of $\mathcal{O}_2$ (orange dot). (b) Frequency $\nu$ and (c) eigenvalues $\Lambda$ associated with the unstable direction of the family of invariant curves as a function of the minimum distance between the fixed point $\mathbf{z}^{\star}$ and the invariant curve $\mathbf{x}(\theta)$. The orange dots are the frequency and the eigenvalue of the periodic orbit $\mathcal{O}_2$.  The parameters are $E_0=0.0925$ and $\omega = 0.0584$.}
\label{fig:nu_lambda}
\end{figure}
\begin{figure}
\centering
\includegraphics[width=0.7\textwidth]{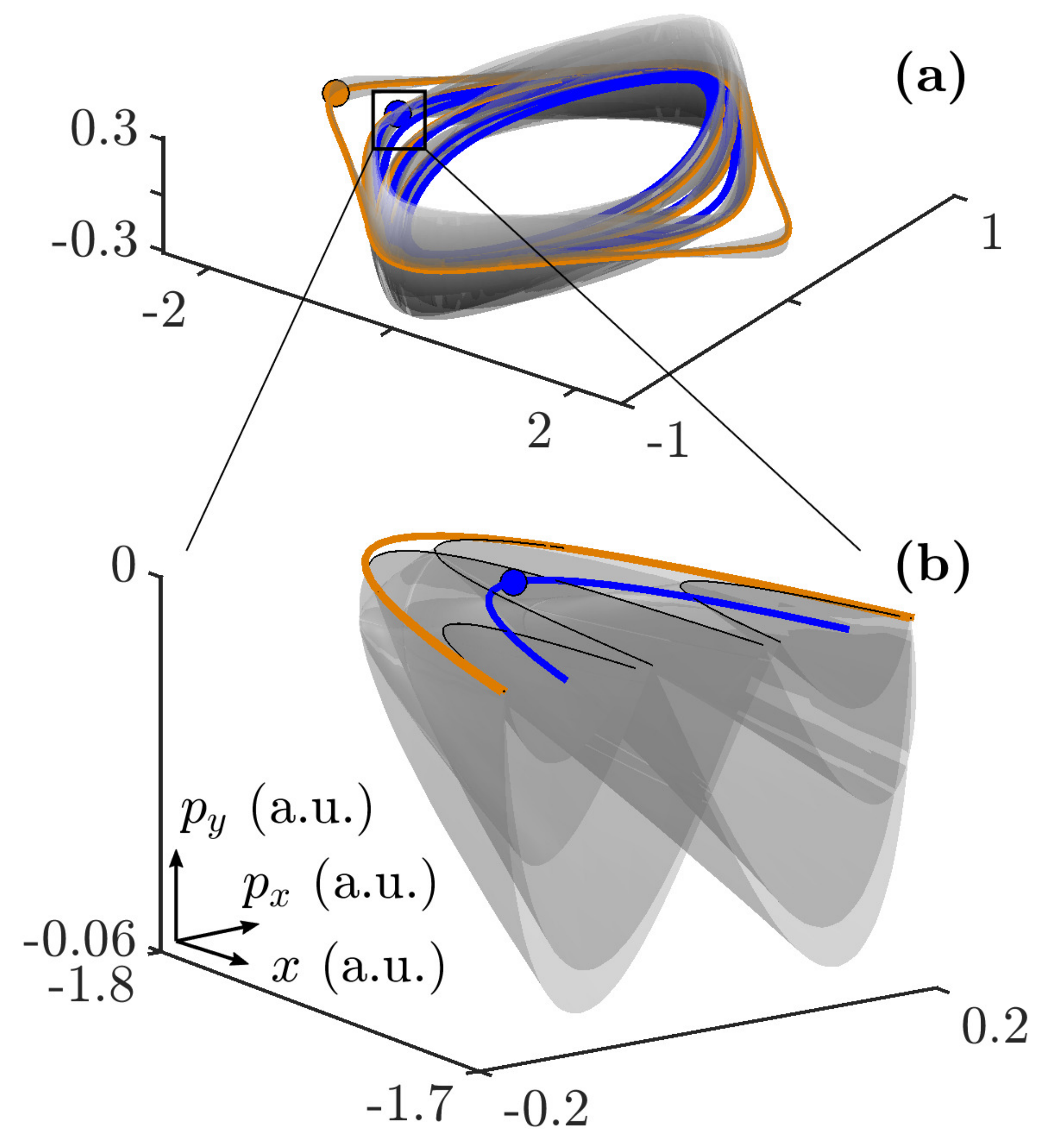}
\caption{Lowest frequency invariant torus of the family ${\cal F}$ for Hamiltonian~\eqref{eq:Hamiltonian_2D}. The orange and blue lines are the periodic orbits $\mathcal{O}_2$ of resonance 7 and 9, respectively. The orange and blue dots are their fixed point under $\mathcal{P}$. The panel (b) is a zoom of (a). In (b), the thin lines are the invariant torus for $p_y = 0$ [note it does not correspond to the intersection with the invariant subspace due to the projection onto the space $(x , p_x , p_y)$]. The parameters are $E_0=0.0925$ and $\omega = 0.0584$.}
\label{fig:tori}
\end{figure}

Natural candidates are the manifolds of the periodic orbits identified in the one-dimensional 
case.
We first examine the linear stability of $\mathcal{O}_1$, $\mathcal{O}_1^{\pm}$ and $\mathcal{O}_2$ (see \cref{tab:eigenvalues_PO1D}). The noticeable feature is that the periodic orbits $\mathcal{O}_1$ and $\mathcal{O}_1^{\pm}$ are hyperbolic-hyperbolic, but ceases to be weakly hyperbolic since they all exhibit an eigenvalue of modulus much greater than one in the transverse direction. 
For instance, the periodic orbit $\mathcal{O}_1$ is about hundreds of times more unstable in the transverse direction than along the polarization direction for this set of parameters.  An electron initiated close to the invariant subspace $y=p_y=0$ and to the unstable manifold of $\mathcal{O}_1$ is immediately pushed away from the ionic core. Given the large excursions of an electron along $\mathcal{O}_1$, the coupling between the degrees of freedom $(x,p_x)$ and $(y,p_y)$ is weak at instances along the periodic orbit, and therefore nearby trajectories tend to escape from the orbit in a nearly-free motion as in the absence of Coulomb potential [i.e., as in Eq.~\eqref{eq:transverse_motion_SFA}]. 
A return to the ionic core within a few tens of laser cycles, relevant for attosecond science~\cite{Corkum2007}, is very unlikely. 
Furthermore, the invariant manifolds associated with $\mathcal{O}_1$ are of codimension two for $\mathcal{P}$. This manifold intersects the plane $y = p_x = 0$ in a zero-dimensional set, e.g., in a discrete set of points if any. Its dimension is too low for this structure to be relevant for the layers observed in \cref{fig:x_py_surface}. The same conclusions can be drawn for $\mathcal{O}_1^{\pm}$. As a consequence, by considering the transverse direction, we can rule out the importance of the periodic orbits $\mathcal{O}_1$ and $\mathcal{O}_1^{\pm}$ in organizing the recollision dynamics. What made $\mathcal{O}_1$ and $\mathcal{O}_1^{\pm}$ so noticeable in the one-dimensional dynamics, namely their large excursions away from the ionic core to reproduce recollision features, became their main weakness in higher dimension (since the large excursions allowed for a fast escape in the transverse directions). RPOs are not necessarily restricted to the polarization plane: There exist some off-polarization axis RPOs for $d=2$. Examples are given in \cref{fig:RPO_2D}. These RPOs are less hyperbolic than the RPOs in the polarization plane (see \cref{tab:eigenvalues_PO2D}) despite the presence of large excursions away form the ionic core, but they remain hyperbolic-hyperbolic, so not driving recollisions in 2D.  

\begin{figure}
\centering
\includegraphics[width=0.8\textwidth]{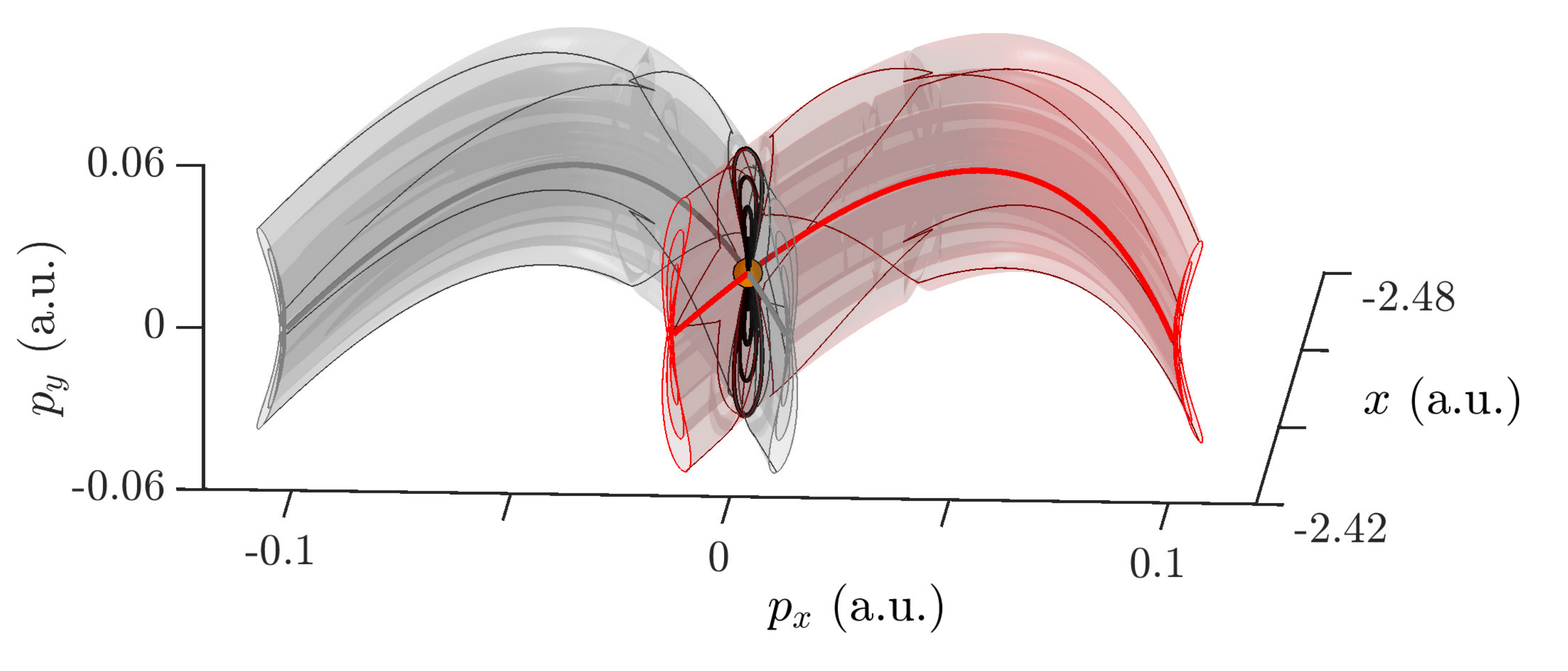}
\caption{Projection onto the space $(x,p_x,p_y)$ of the stable (red) and unstable (gray) manifolds of the family of invariant tori $\mathcal{F}$ (black curves) associated with the hyperbolic-elliptic periodic orbit $\mathcal{O}_2$ (orange marker) under the map $\mathcal{P}$ for Hamiltonian~\eqref{eq:Hamiltonian_2D}. The dark gray and dark red curves along the outer branches of the manifolds are a set of trajectories on the manifold initiated close to the invariant curves. The red and gray thick curves are the stable and unstable manifolds of the hyperbolic-elliptic periodic orbit $\mathcal{O}_2$ also shown on \cref{fig:distance_manifold}c. The parameters are $E_0=0.0925$ and $\omega = 0.0584$. }
\label{fig:O2_manifolds}
\end{figure}
The second candidate, the periodic orbit $\mathcal{O}_2$ is weakly hyperbolic in the polarization direction and elliptic in the transverse direction (see \cref{tab:eigenvalues_PO1D}). 
In the elliptic direction, 
$\mathcal{O}_2$ gives rise to a 
Cantor family $\mathcal{F}$ of hyperbolic invariant tori \cite{Graff74,BroerHS96}. 
Given that Hamiltonian~\eqref{eq:Hamiltonian_2D} is analytic, the size of the holes of the Cantor structure 
are exponentially small with the order of the resonance \cite{JorbaV97b,JorbaV97a}. As a consequence, only low order resonances produce holes that are visible under 
the usual computer resolution. In our case there are no 
resonances of low order near $\mathcal{O}_2$ and, 
hence, the Cantor family $\mathcal{F}$ of tori looks like a continuous family.
The invariant curves are computed using the method described in Ref.~\cite{Jorba2001} and 
summarized in Appendix A. The black curves are invariant curves surrounding $\mathcal{O}_2$. A few of these invariant curves are depicted in \cref{fig:nu_lambda}. They correspond to two-dimensional tori for the continuous Hamiltonian flow, and one such two-dimensional torus is represented in \cref{fig:tori} together with the relevant SPOs. We notice that the family ${\cal F}$ of invariant tori wraps around ${\cal O}_2$. The frequency of these invariant curves as a function of the distance between ${\bf z}^\star$ and the invariant curve is represented in \cref{fig:nu_lambda}b, where ${\bf z}^\star$ is the fixed point of $\cal P$ associated with $\mathcal{O}_2$ (orange dot in the figures). This family of invariant curves seems to end when a major resonance is encountered, here $1:9$ for the chosen parameters (blue dashed line in \cref{fig:nu_lambda}b). The SPO associated with the $1:9$ resonance is depicted in blue in \cref{fig:tori}.  
In \cref{fig:nu_lambda}c, we plotted the eigenvalue of each torus as a function between ${\bf z}^\star$ and the invariant curve. We notice that all of these eigenvalues are of order 1, hence this family of invariant tori is weakly hyperbolic. 
Given that $\mathcal{O}_2$ is unstable along the polarization direction, each invariant curve of $\mathcal{F}$ has a stable and an unstable manifold. Pieces of the stable and unstable manifolds of $\mathcal{F}$ in the vicinity of the invariant curves are depicted in \cref{fig:O2_manifolds}. Because the holes of the Cantor family $\mathcal{F}$ of invariant curves are not visible at this scale, $\mathcal{F}$ is seemingly two-dimensional under the map $\mathcal{P}$. Therefore, the invariant manifolds of $\mathcal{F}$ are seemingly three-dimensional, hence of codimension one in phase space. The intersections between the stable manifold of $\mathcal{F}$ and the plane $(x,p_y)$ for $y = p_{x}=0$ are lines. As a consequence, the resulting union of the stable manifolds of the tori of $\cal F$, denoted $\mathcal{W}_{\cal F}^{\rm s}$ in what follows, has the desired dimension. 

In order to show that the invariant manifolds of the family of invariant curves $\mathcal{F}$ associated with $\mathcal{O}_2$ drives the electron dynamics and the recollisions, we consider their intersections with a slice of phase space. For $d=2$, we resort to slices in the four-dimensional space of initial conditions, in order to represent the distance as a function of the initial conditions. In \cref{fig:x_py_surface} we chose $y=p_x = 0$. It then amounts to computing the intersections between $\mathcal{W}_{\cal F}^{\rm s}$ and this slice of phase space of codimension $2$. Note that this is, in principle, possible since $\mathcal{W}_{\cal F}^{\rm s}$ is of dimension $3$. Below, we detail the procedure we follow to determine these intersections.

\subsubsection{Computation of the intersection of $\mathcal{W}_{\cal F}^{\rm s}$ and a slice of phase space}

Each invariant curve of the family $\cal F$ is labeled by $\sigma$ and its representation is denoted $\mathbf{x}_{\sigma}$. We focus on the intersections between $\mathcal{W}_{\cal F}^{\rm s}$ and the slice $y=p_x=0$ for $d=2$. First, we consider the stable manifold of one invariant curve $\mathbf{x}_{\sigma}$ with frequency $\nu$. The real eigenvalue related to the stable direction of the invariant curve is denoted $\Lambda^{\rm s}$ and $\boldsymbol{\Psi}_{\sigma}^{\rm s}(\theta)$ is the vector tangent to the stable manifold at each point $\mathbf{x}_{\sigma} (\theta)$ of the invariant curve (see Ref.~\cite{Jorba2001} and \cref{sec:linear_stability_curves}). 
The linear approximation to the stable manifold is given by
\begin{equation}
\label{eq:linear_approximation_invariant_manifolds_invariant_tori}
\mathbf{z}_{\sigma} (s,\theta) = \mathbf{x}_{\sigma} (\theta) + s \boldsymbol{\Psi}_{\sigma}^{\rm s} (\theta) ,
\end{equation}
where $s$ is a real parameter along the manifold and $\theta \in \mathbb{T}=\mathbb{R} / 2\pi \mathbb{Z}$.
A fundamental domain of the invariant manifold related to an invariant curve 
is a minimal cylinder from which the entire stable manifold can 
be obtained by iterating the map ${\cal P}^{-1}$. For instance, in the linear approximation given by Eq.~\eqref{eq:linear_approximation_invariant_manifolds_invariant_tori}, a fundamental domain is the set of points $\mathbf{z}_{\sigma} (s, \theta)$ such that $(s,\theta) \in \{ [h, h/\Lambda^{\rm s}]\times \mathbb{T}\}$, where $h\in \mathbb{R}$. 
In order to produce a good enough approximation to a fundamental domain, we select a sufficiently small value of $h$ such that the linear approximation satisfies the invariance equation $| {\cal P}^{-1} (\mathbf{z}_{\sigma} (s,\theta)) - \mathbf{z}_{\sigma} (s/\Lambda^{\rm s},\theta - \nu) | < \epsilon $, where $\epsilon$ is a suitable threshold. 
As $\mathcal{W}^{\rm s}_{\cal F}$ denotes the union of all the stable manifolds of the invariant curves 
of the family $\cal F$, a fundamental domain of $\mathcal{W}^{\rm s}_{\cal F}$ is given by the union of the fundamental domains associated with each manifold, and can be parametrized by $\sigma$, $s$ and $\theta$. This makes it an object of dimension $3$. 

The stable manifold related to the invariant curve $\mathbf{x}_{\sigma}$ intersects the plane 
$y=p_x=0$ in one or multiple points for $s = s^{\star}$ and $\theta = \theta^{\star}$. This condition reads 
\begin{subequations}
\label{eq:Newton_method_intersection_unstable_plane_2D}
\begin{eqnarray}
{\cal P}^{-m} (\mathbf{z}_{\sigma} (s^{\star},\theta^{\star})) \cdot \mathbf{e}_y = 0 , \\
{\cal P}^{-m} (\mathbf{z}_{\sigma} (s^{\star},\theta^{\star})) \cdot \mathbf{e}_{p_x} = 0 .
\end{eqnarray}
\end{subequations}
The vectors $\mathbf{e}_y$ and $\mathbf{e}_{p_{x}}$ are unit vectors along the direction $y$ and $p_{x}$, 
respectively. The parameter $m$ corresponds to the number of iterations required for the trajectory to 
reach the plane $y = p_x = 0$. The event that the trajectories from the fundamental domain of 
$\mathcal{W}^{\rm s}_{\cal F}$ is in the neighborhood of the plane is rare because 
the plane $y=p_x=0$ is codimension $2$. In order to remedy this difficulty, we proceed as follow:
\begin{enumerate}
\item
We compute a finite set of invariant curves $\mathbf{x}_{\sigma}$ associated with $\mathcal{O}_2$ such that 
$\sigma = 1, \hdots , N_{\sigma}$. We have used $N_{\sigma} \sim 100$. Then, we compute, for each 
$\theta$, the stable direction $\boldsymbol{\Psi}^{\rm s}_{\sigma} (\theta)$ associated with 
the real eigenvalue of modulus smaller than unity.
\item
For each invariant curve, i.e., for each $\sigma$, we consider a mesh of initial conditions in the fundamental domain denoted $\mathbf{z}_{\sigma} (s_i , \theta_j)$ given by Eq.~\eqref{eq:linear_approximation_invariant_manifolds_invariant_tori} for all $\sigma$, and for $i = 0, \hdots , N_s$ and $j = 0 , \hdots , N_{\theta}$. We obtain an ensemble of initial conditions $\mathbf{z}_{\sigma} (s_i,\theta_j)$ labeled by $\sigma$, $i$ and $j$. We have used $N_s \sim 100$ and $N_{\theta} \sim 100$. As a consequence, there are approximately $10^6$ initial conditions in the fundamental domain of $\mathcal{W}_{\cal F}^{\rm s}$. 
\item
We compute ${\cal P}^{-m} ( \mathbf{z}_{\sigma} ( s_i ,\theta_j ) )$ for $m = 1, \hdots , N_m$. We have used $N_m \sim 100$. We store the points which are close to the plane $y = p_x = 0$ and their associated $m$, $s_i$, $\theta_j$ and $\sigma$. A point is stored if it is close enough to the plane $y = p_x = 0$. We have used a tolerance of $10^{-2}$. This tolerance is arbitrary and mainly influences the acceptation/rejection rate in the next step.
\item
The points which are stored in the previous step are refined using a Newton method. We use a Newton method to compute $s^{\star}$ and $\theta^{\star}$ inside the fundamental domain such that for a given $m$ and $\sigma$, Eq.~\eqref{eq:Newton_method_intersection_unstable_plane_2D} is satisfied.
The two parameters to adjust in order to fulfill this condition are $s$ and $\theta$.
The points and their labels which have been stored in the previous step are initial guesses for the Newton method.
Note that at each iteration of the Newton method, the derivatives of Eqs.~\eqref{eq:Newton_method_intersection_unstable_plane_2D} with respect to $s$ and $\theta$ must be evaluated and must be preliminary derived from the explicit definition of $\mathbf{z}_{\sigma} (s,\theta)$ given by Eq.~\eqref{eq:linear_approximation_invariant_manifolds_invariant_tori}. As a result, the images of the points initiated on the stable manifold, which belong to the stable manifold, also belong to the plane $y = p_x = 0$. Therefore these points correspond to the intersection between $\mathcal{W}_{\cal F}^{\rm s}$ and the slice $y=p_x=0$. The points such that $s > h/ \Lambda^{\rm s}$ are not kept since they are outside the fundamental domain of $\mathcal{W}_{\cal F}^{\rm s}$.
\end{enumerate}
The same scheme can be used to compute the intersection between the plane $y = p_x = 0$ with the unstable manifold $\mathcal{W}^{\rm u}_{\cal F}$ by changing $\Lambda^{\rm s} \to 1/\Lambda^{\rm u}$ (which are actually the same values, i.e., $\Lambda^{\rm s} = 1/\Lambda^{\rm u}$), $\boldsymbol{\Psi}_{\sigma}^{\rm s} \to \boldsymbol{\Psi}_{\sigma}^{\rm u}$, $\nu \to - \nu$ and ${\cal P}^{-1} \to {\cal P}$. 

\subsubsection{Results: $\mathcal{F}$ drives the recollisions}
In \cref{fig:x_py_surface}, the red dots are the numerically determined points on the intersection between the plane of initial conditions and the stable manifold of $\mathcal{F}$.  We see that the red dots are aligned along the chaotic layers shown on the contour plot of the final distance of the electron. It is a strong evidence that the stable manifold of $\mathcal{F}$ structures the recollision dynamics: The return of the electron is guided by the stable manifold of $\mathcal{F}$ that gets close to the ionic core, and its escape is through its unstable manifold. Also, $\mathcal{F}$ organizes the electron dynamics by separating the phase space into regions of bounded and unbounded motion, recolliding or directly ionizing electrons. The manifolds of $\cal F$ draw codimension-one ``tubes'' around the subspace $y=p_y=0$ ensuring that the recolliding trajectories remain along the polarization plane, and as a consequence, that the results obtained for one dimension are structurally stable. 
\par
Close to the intersection of the stable and the unstable manifolds of $\mathcal{F}$, there is an infinite number of periodic orbits with large period due to the back and forth motion along the invariant manifolds. These orbits are RPOs as we can see in the 1D case in \cref{fig:distance_manifold}. The recollision scenario in 2D is fully compatible with the scenario in 1D. The weak instability of the family of invariant tori ensures that the recolliding trajectories spend a significant amount of time along the polarization axis. This reinforces the relevance of the reduction to the dynamics along the polarization axis. 

\subsection{Recollisions in 3D}

For $d=3$, Hamiltonian~\eqref{eq:Hamiltonian} written in polar coordinates along the $\hat{\mathbf{x}}$-direction becomes
\begin{equation}
\label{eq:Hamitlonian_3D}
H (x , \rho , \theta , p_x , p_{\rho} , p_{\theta},t ) = \dfrac{p_x^2}{2} + \dfrac{p_{\rho}^2}{2} + \dfrac{p_{\theta}^2}{2 \rho^2} - \dfrac{1}{\sqrt{x^2 + \rho^2 + 1}} + x E_0 \cos (\omega t) . 
\end{equation}
This set of coordinates highlights a conserved quantity $p_{\theta}=y p_z-z p_y$, the angular momentum in the plane transverse to the polarization axis. Therefore, the phase space is foliated by leaves with constant values of $p_{\theta}$, and periodic orbits are parabolic in the transverse direction (which explains the eigenvalue 1 in \cref{tab:eigenvalues_PO1D} and \cref{tab:eigenvalues_PO2D}). As a consequence, a small perturbation from the initial conditions along this direction does not modify drastically the electron dynamics. Therefore, the recollision scenario described for $d=2$ remains identical for $d=3$.

\section{Conclusions}

In summary, the stable and unstable manifolds of the family of invariant tori $\mathcal{F}$ of $\mathcal{O}_2$, structure the recollisions and organize the electron dynamics in phase space. On one side, the stable and unstable manifolds surround the region of bounded motion, and their intersection produces delays in the ionization of the electron known as RESI~\cite{Mauger2012_PRL,Mauger2012_PRE}. On the other side, the stable and unstable manifolds are the skeleton of the recollision dynamics by channeling the electrons first away from the ionic core and then back, following the stable manifold $\mathcal{W}_{\cal F}^{\rm s}$. The ``tubes'' drawn by the manifolds of $\cal F$ around the invariant subspace $y=p_y=0$ ensure the validity of the results obtained in the one-dimensional case along the polarization axis.   

\appendix

\section{Computation of the family $\mathcal{F}$ of invariant curves and their linear stability}

In this section, we briefly summarize a method to compute invariant curves and their linear stability as explained in Ref.~\cite{Jorba2001}. For other methods, we refer to Refs.~\cite{DieciLR91, Simo98, SchilderOV05, JorbaO09, HaroCLMF16}.
For simplicity, we denote the phase-space variables $\mathbf{z} = (\mathbf{r} , \mathbf{p})$ and $\mathbf{z} \in \mathbb{R}^N$ with $N = 2d$. We define the Poincar\'e / stroboscopic map ${\cal P}$ such that
\begin{equation}
\label{eq:stroboscopic_map}
{\cal P}( \mathbf{z}) = {\phi}_{0}^{T} (\mathbf{z} ) ,
\end{equation}
where ${\phi}_0^T$ is the Hamiltonian flow over one period of the laser field. We reduce the study of the continuous Hamiltonian flow to the study of a map
\begin{equation}
\label{eq:reduced_dynamical_system}
\bar{\mathbf{z}} = {\cal P} (\mathbf{z}) .
\end{equation}
The periodic orbits of Hamiltonian~\eqref{eq:Hamiltonian} are fixed points of the map~\eqref{eq:reduced_dynamical_system}, and the two-dimensional invariant tori are invariant curves of ${\cal P}$. 

We denote $\mathbf{x}$ a closed curve parametrized by $\theta \in \mathbb{T}=\mathbb{R} / 2\pi \mathbb{Z}$. If the curve is invariant by ${\cal P}$, there exists $\nu \in \mathbb{R} \setminus 2 \pi \mathbb{Q}$ such that
\begin{equation}
\label{eq:zero_invariant_curve}
{\cal P}(\mathbf{x}(\theta)) = \mathbf{x} (\theta + \nu) .
\end{equation}
We use a Fourier representation method~\cite{Jorba2001} to compute the invariant curves of the map ${\cal P}$ and their linear stability. The Fourier representation of the invariant curve is
\begin{equation}
\label{eq:fourier_representation_invariant_curve}
\mathbf{x}(\theta ; \boldsymbol{\alpha} ) = \mathbf{a}_0 + \sum_{k=1}^M \left[ \mathbf{a}_k \cos (k \theta) + \mathbf{b}_k \sin (k\theta) \right] ,
\end{equation}
where $2 M + 1$ is the number of coefficients and we denote $\boldsymbol{\alpha} = ( \mathbf{a}_0, \hdots , \mathbf{a}_M , \mathbf{b}_1 , \hdots , \mathbf{b}_M )$. For numerical purposes, $M$ is large enough so that $|\mathbf{a}_M| + |\mathbf{b}_M|$ is small. 
Computing invariant curves consists in determining the set of Fourier coefficients $\boldsymbol{\alpha}$ for which Eq.~\eqref{eq:zero_invariant_curve} is satisfied.

The Poincar\'{e} map for the Hamiltonians we consider does not explicitly depend on the parameter of the invariant curve $\theta$. This case is referred to as the \emph{autonomous} case~\cite{Jorba2001}. In this case, an infinity of representations of the invariant curves is possible, and an additional condition must be imposed in order to fix the representation.
For pedagogical reasons, the method is first described for the non-autonomous case, for which ${\cal P}$ depends explicitly on the parameter $\theta$, then for the autonomous case, for which  ${\cal P}$ does not explicitly depend on $\theta$.
Then, we summarize the method to compute the linear stability of the invariant curves.

\subsection{Non-autonomous maps}

First, we consider that ${\cal P}$ depends explicitly on $\theta$. The dynamical system~\eqref{eq:reduced_dynamical_system} is written as a skew-product 
\begin{eqnarray*}
\overline{\mathbf{z}} &= {\cal P} (\mathbf{z}, \theta) , \\
\overline{\theta} &= \theta + \nu .
\end{eqnarray*}
In Eq.~\eqref{eq:zero_invariant_curve}, we use the Fourier representation of the invariant curve given by Eq.~\eqref{eq:fourier_representation_invariant_curve}. There are $2 M + 1$ coefficients, and as a consequence, $(2 M +1) N$ unknowns. We consider a mesh of parameters $\theta = \theta_0 , \hdots , \theta_{2M}$ such that
\begin{equation}
\label{eq:discretization_theta_invariant_curves_computation}
\theta_j = j \frac{2 \pi}{2 M+1} ,
\end{equation}
with $0 \leq j \leq 2 M$. Each point $\mathbf{x} (\theta_j ; \boldsymbol{\alpha})$ is a point on the curve, so it is a solution of Eq.~\eqref{eq:zero_invariant_curve}. Hence, the zero function of the $2 M+1$ points solution of the invariance equation~\eqref{eq:zero_invariant_curve} reads explicitly 
\begin{equation}
\label{eq:invariant_curves_zero_function_non_autonomous}
{F}_{\rm na}(\boldsymbol{\alpha} , \nu) = 
\begin{bmatrix}
{\cal P} \left( \mathbf{x} (\theta_0 ; \boldsymbol{\alpha}) \right) - \mathbf{x} (\theta_0 + \nu ; \boldsymbol{\alpha})  \\
\vdots \\
{\cal P} \left( \mathbf{x} (\theta_{2M} ; \boldsymbol{\alpha}) \right) - \mathbf{x} (\theta_{2M} + \nu ; \boldsymbol{\alpha}) 
\end{bmatrix} = \boldsymbol{0} .
\end{equation}
Equation~\eqref{eq:invariant_curves_zero_function_non_autonomous} is a set of $(2M+1)N$ equations where the Fourier coefficients $\boldsymbol{\alpha}$ are unknowns.
A Newton method is used to find numerically the coefficients $\boldsymbol{\alpha}$ solution of Eq.~\eqref{eq:invariant_curves_zero_function_non_autonomous}.

We denote $\mathbf{z}^{\star}$ the location of the fixed point under ${\cal P}$.
The initial guess of the Newton method is given by $\mathbf{a}_0 = \mathbf{z}^{\star}$ and $\mathbf{a}_j = \mathbf{b}_j = \boldsymbol{0}$ for $j = 2 , \hdots , M$. The first coefficients are $\mathbf{a}_1 = h \Re ( \mathbf{v}_{\rm c} )$ and $\mathbf{b}_1 = h \Im ( \mathbf{v}_{\rm c} )$, or $\mathbf{a}_1 = h \Im ( \mathbf{v}_{\rm c} )$ and $\mathbf{b}_1 = h \Re ( \mathbf{v}_{\rm c} )$, depending on the running direction of the electrons along the invariant curve. The complex vector $\mathbf{v}_c$ is the eigenvector associated with the center component of the fixed point (i.e., associated with the complex eigenvalue of the fixed point), see also \cite{GabernJ04a}. 

\subsection{Autonomous maps and continuation method}

When the map ${\cal P}$ does not depend explicitly on the parameter $\theta$, an additional condition must be imposed to compute the Fourier coefficients of the invariant curve for its Fourier representation is not unique. Indeed, let $\mathbf{x} (\theta ; \boldsymbol{\alpha} ) = \mathbf{y}(\theta + \phi ; \boldsymbol{\alpha} )$ with $\phi \in \mathbb{T}$. The invariant curves ${\bf x}$ and ${\bf y}$ parametrized by $\theta$ have different Fourier series but represent the same curve, the origin $\theta = 0$ of the parametrization being shifted in the $\theta$-space. There are many different ways to overcome this difficulty. Here we choose to fix arbitrarily one component of the origin of the invariant curve $\mathbf{x} (0 ; \boldsymbol{\alpha}) = \sum_{k=0}^M \mathbf{a}_k$. Given a unit vector $\bf e$, we consider the additional condition
\begin{equation}
\label{eq:autonomous_extra_equation}
F_{\rm e} (\boldsymbol{\alpha}) = \mathbf{e} \cdot \left[ \mathbf{x} (0 ; \boldsymbol{\alpha}) - \mathbf{z}^{\star} \right] = 0 .
\end{equation}
Choosing a relevant vector $\mathbf{e}$ depends on the system we consider. In our case, we fix the origin of the invariant curve along $p_x$, i.e., for $d=2$, we choose $\mathbf{e} = (0,0,1,0)$.
In the autonomous case, the coefficients of the Fourier series of the invariant curves are also zeros of the function $F_{\rm e}$. The Fourier coefficients are solution of Eqs.~\eqref{eq:invariant_curves_zero_function_non_autonomous} and \eqref{eq:autonomous_extra_equation}. They are determined using a Newton method. 
As a consequence of this extra condition, the linear system we solve in the Newton method becomes non-square, i.e., there are more conditions than unknowns. In order to solve the non-square linear system, we use a Gaussian elimination with maximal pivoting. In this way, the redundant equation is sent to the last row, and reads ``$0 = 0$''.

In addition, in the autonomous case, invariant curves come in families, labeled by different irrational frequencies $\nu$.
A continuation method is used to compute the invariant curves of the family. 
The frequency of each invariant curve $\nu$ is a free parameter in the continuation method, and as a consequence the unknowns are $\boldsymbol{\alpha}$ and $\nu$. Note that $\nu$ enters into play in the function $F_{\rm na}$ in the non-autonomous case and one must differentiate this function with respect to $\nu$ in the Newton method.

The continuation method consists in computing the invariant curves one after the others. A parameter $\delta$ controls the distance between two invariant curves. We distinguish two cases: ($i$) no invariant curves are computed yet, only the fixed point is known, and ($ii$) at least one invariant curve of the family is known. In the latter case, the Fourier coefficients of the last computed invariant curve are denoted $\tilde{\boldsymbol{\alpha}}$. The coefficients of the new invariant curve $\boldsymbol{\alpha}$ are such that
\begin{equation}
\label{eq:continuation_method_invariant_curves_autonomous}
F_{\delta} (\boldsymbol{\alpha} ) = 
\left\lbrace 
\begin{array}{l @{\qquad} l}
| \mathbf{a}_1 |^2 + | \mathbf{b}_1 |^2 - \delta^2 = 0 , & \text{if } (i) , \\
| \boldsymbol{\alpha} - \tilde{\boldsymbol{\alpha}} |^2 - \delta^2 = 0 , & \text{if } (ii) .
\end{array}
\right.
\end{equation}
The above-condition ($i$) ensures the curves to converge sufficiently far away from the fixed point. 
In the Newton method, the initial guess is determined from interpolating the Fourier coefficients and frequencies $\nu$ of the previously computed invariant curves of the family $\mathcal{F}$. In our computations, we use a linear interpolation.

In sum, in the autonomous case, the family of invariant curves associated with a fixed point $\mathbf{z}^{\star}$ is computed by solving numerically the system of $N (2 M + 1) + 2$ equations [see Eqs.~\eqref{eq:invariant_curves_zero_function_non_autonomous}, \eqref{eq:autonomous_extra_equation} and \eqref{eq:continuation_method_invariant_curves_autonomous}] and $N (2 M +1) + 1$ unknowns $\boldsymbol{\alpha}$ and $\nu$. Among the $N (2 M + 1) + 2$ equations, $N (2 M+1)$ equations correspond to the non-autonomous case~\eqref{eq:invariant_curves_zero_function_non_autonomous} with the function ${F}_{\rm na} (\boldsymbol{\alpha} , \nu)$. One extra equation is used to fix the non-uniqueness of the Fourier representation of the invariant curve. It is given by Eq.~\eqref{eq:autonomous_extra_equation} and the function $F_{\rm e} (\boldsymbol{\alpha})$. The frequency of the invariant curve $\nu$ is a free parameter, and $\delta$ is used to maneuver along the family of invariant tori. It is given by Eq.~\eqref{eq:continuation_method_invariant_curves_autonomous} and the function $F_{\rm \delta} (\boldsymbol{\alpha})$. 

\subsection{Linear stability of the invariant curves \label{sec:linear_stability_curves}}

We consider a given invariant curve $\mathbf{x}$ determined by the method summarized above.
The method for computing the linear stability of an invariant curve is the same for autonomous and non-autonomous maps. We denote $\boldsymbol{\Psi} (\theta)$ the eigenvectors along the invariant curve and $\Lambda$ its associated eigenvalue. For each $\theta$, we consider a small displacement from the invariant curve $\delta \mathbf{x} (\theta)$. The invariance equation~\eqref{eq:zero_invariant_curve} becomes ${\cal P} (\mathbf{x} (\theta) + \delta \mathbf{x} (\theta)) = \mathbf{x} (\theta+\nu) + \delta \mathbf{x} (\theta+\nu)$. We perform a Taylor expansion on the left-hand side, we obtain ${\cal P} (\mathbf{x} (\theta) + \delta \mathbf{x} (\theta)) \approx {\cal P} (\mathbf{x} (\theta)) + \mathbb{B} (\theta) \delta \mathbf{x} (\theta)$, where $\mathbb{B} (\theta) = \partial {\cal P} (\mathbf{x} (\theta)) / \partial \mathbf{z}$ which corresponds to the tangent flow at each point of the invariant curve. The eigenvalue problem to solve is 
\begin{equation}
\label{eq:CHAP5_stability_invariant_curve}
\mathbb{B} (\theta) \boldsymbol{\Psi} (\theta) = \Lambda \boldsymbol{\Psi} (\theta + \nu ) .
\end{equation}
\par
We consider the discrete value problem for $2 M+1$ angles $\theta_j$ defined in Eq.~\eqref{eq:discretization_theta_invariant_curves_computation}. The Fourier representation of the eigenvectors is denoted
\begin{equation*}
\boldsymbol{\Psi} (\theta ; \underline{\boldsymbol{\alpha}}) = \mathbf{A}_0 + \sum_{k=1}^M \left[ \mathbf{A}_k \cos (k \theta) + \mathbf{B}_k \sin (k \theta) \right] ,
\end{equation*}
where the Fourier coefficients of the eigenvector are $\underline{\boldsymbol{\alpha}} = (\mathbf{A}_0 , \hdots , \mathbf{A}_M , \mathbf{B}_1 , \hdots , \mathbf{B}_M )$. The rotation matrix of an angle $\nu$ is denoted $\mathbb{H}_{\nu}$ and is such that $\mathbb{H}_{\nu} \boldsymbol{\Psi} (\theta ; \underline{\boldsymbol{\alpha}}) = \boldsymbol{\Psi} (\theta + \nu ; \underline{\boldsymbol{\alpha}} )$.
The matrix $\mathbb{H}_{\nu}$ is of dimension $[(2 M+1) N] \times [(2 M +1) N]$. Note that, due to the Fourier representation of the invariant curve, the matrix $\mathbb{B}(\theta) = \partial \boldsymbol{\cal P} (\mathbf{x} (\theta)) / \partial \mathbf{z}$ has to be determined using the chain rule.
Finally, in Fourier representation, the eigenvalues $\Lambda$ and eigenvectors  $\boldsymbol{\Psi} ( \theta ; \underline{\boldsymbol{\alpha}} )$ of the invariant curve $\mathbf{x}(\theta;\boldsymbol{\alpha})$ are solutions of the eigenvalue problem 
\begin{equation*}
\mathbb{H}_{\nu}^{\dagger} \mathbb{B} \boldsymbol{\Psi} ( \theta ; \underline{\boldsymbol{\alpha}} ) = \Lambda  \boldsymbol{\Psi} ( \theta ; \underline{\boldsymbol{\alpha}} ) .
\end{equation*}
As a result, there are $(2 M + 1)N$ eigenvalues. In our case, as $N=4$, ${\cal P}$ is autonomous and each torus is hyperbolic, four of these eigenvalues are real. Two of them are equal to the unity and correspond to 
the tangent direction of the invariant curve and the other one corresponds the direction of the family.
The other two real eigenvalues correspond to the eigenvalues of the invariant curve. 
One is of modulus smaller than one and corresponds to the eigenvalue associated with the stable direction 
of the invariant curve, and the other one has a modulus greater than one corresponding 
to the unstable direction.

\section*{Acknowledgments}
The project leading to this research has received funding from the European Union’s Horizon 2020 research and innovation program under the Marie Sk\l odowska-Curie Grant Agreement No. 734557. MJC and AJ have been supported by the Spanish grants
PGC2018-100699-B-I00 (MCIU/AEI/FEDER, UE) and the
Catalan grant 2017 SGR 1374.


\newcommand{\noopsort}[1]{}

\end{document}